\documentstyle[11pt]{article}
\begin{document}

\title{{\bf General Kastler-Kalau-Walze Type Theorems for Manifolds with  Boundary II}
\thanks{partially supported by NSFC 11771070.
 }}
\author{ Yong Wang\\
{\scriptsize \it School of Mathematics and Statistics, Northeast Normal University, Changchun Jilin 130024,  China ;}\\
{\scriptsize \it E-mail: wangy581@nenu.edu.cn}}

\date{}
\maketitle

\noindent {\bf Abstract}~~In this paper, we establish some general Kastler-Kalau-Walze type theorems
for any dimensional manifolds with  boundary which generalize the results in [WW1].
\\
\noindent{\bf MSC:}\quad 58G20; 53A30; 46L87\\
 \noindent{\bf Keywords:}~
Noncommutative residue for manifolds with boundary; Dirac operators\\
\section{Introduction}
\quad The noncommutative residue found in [Gu],[Wo] plays a prominent role in noncommutative geometry.
For one-dimensional manifolds, the noncommutative residue was discovered by Adler [Ad]
 in connection with geometric aspects of nonlinear partial differential equations.
 For arbitrary closed compact $n$-dimensional manifolds, the noncommutative residue was introduced by Wodzicki in [Wo] using the theory of zeta
 functions of elliptic pseudodifferential operators.
In [Co1], Connes used the noncommutative residue to derive a conformal 4-dimensional Polyakov action analogue.
Furthermore, Connes made a challenging observation that the noncommutative residue of the square of the inverse of the
Dirac operator was proportional to the Einstein-Hilbert action in [Co2].
Let $s$ be the scalar curvature and Wres denote  the noncommutative residue. Then the Kastler-Kalau-Walze
theorem gives an operator-theoretic explanation of the gravitational action and says that for a $4-$dimensional closed spin manifold,
 there exists a constant $c_0$, such that
 $$
{\rm  Wres}(D^{-2})=c_0\int_Ms{\rm dvol}_M.\eqno(1.1)
$$
In [Ka], Kastler gave a brute-force proof of this theorem. In [KW], Kalau and Walze proved this theorem in the
normal coordinates system simultaneously. And then, Ackermann proved that
the Wodzicki residue ${\rm  Wres}(D^{-2})$ in turn is essentially the second coefficient
of the heat kernel expansion of $D^{2}$ in [Ac].

On the other hand, Fedosov etc. defined a noncommutative residue on Boutet de Monvel's algebra and proved that it was a
unique continuous trace in [FGLS]. In [Sc], Schrohe gave the relation between the Dixmier trace and the noncommutative residue for
manifolds with boundary.
 We [Wa2] proved a Kastler-Kalau-Walze type theorem for $4$-dimensional spin manifolds with boundary.
  Furthermore, We [Wa3] generalized the definition
of lower  dimensional volumes in [Po] to  manifolds with boundary and found a Kastler-Kalau-Walze type theorem for
$6$-dimensional manifolds with boundary. In [WW2], we established a general Kastler-Kalau-Walze type theorem
for any dimensional manifolds with boundary which generalized the results in [Wa3].
In [WW1],  we computed the lower
dimensional volume ${\rm Vol}_{6}^{(1,3)}$ for $6$-dimensional spin manifolds with boundary and got another Kastler-Kalau-Walze type theorem for
$6$-dimensional manifolds with boundary.  {\bf The motivation of this paper} is  to establish a general Kastler-Kalau-Walze type theorem
for any dimensional manifolds with boundary which extends the theorem in [WW1].
Our main result is as follows.\\

\noindent{\bf Theorem 1.1}
{\it Let M be a $n=2m+2$-dimensional compact spin manifold  with the boundary $\partial M$, then}
\begin{eqnarray*}
 \widetilde{{\rm Wres}}\big[\pi^+D^{-1}\circ \pi^+D^{(-2m+1)}\big]
 &=& \frac{(2-n)(2\pi)^{{n}/2}}{12\Gamma({n}/2)}\int_Ms{\rm dvol}_M
  \nonumber\\
 &&+{\frac{2}{1-n}} L_{0}Vol(S_{n-2})\int_{\partial_{M}}K{\rm dvol}_{\partial_{M}},
~~~~~~~~~~~~~(1.2)
\end{eqnarray*}
{\it where $Vol(S_{n-2})$ is the canonical volume of the $n-2$ dimensional sphere $S_{n-2}$ and $s$ is the scalar curvature on $M$ and $K$ is the extrinsic curvature on $\partial M$ as well as
 $L_0$ is a constant (see (3.42))}.\\

\indent This paper is organized as follows: In Section 2, we recall the definition of the lower dimensional volumes of compact Riemannian manifolds
 with  boundary. In Section 3, we establish some general Kastler-Kalau-Walze type theorems
for any dimensional manifolds with boundary which extend the theorem in [WW1].

\section{Lower-Dimensional Volumes of Spin Manifolds  with  boundary}

\quad In this section we consider an $n=2m+2$-dimensional oriented Riemannian manifold $(M, g^{M})$ with boundary $\partial{ M}$ equipped
with a fixed spin structure. We assume that the metric $g^{M}$ on $M$ has
the following form near the boundary
 $$
 g^{M}=\frac{1}{h(x_{n})}g^{\partial M}+\texttt{d}x _{n}^{2} ,\eqno(2.1)
$$
where $g^{\partial M}$ is the metric on $\partial M$. Let $U\subset
M$ be a collar neighborhood of $\partial M$ which is diffeomorphic $\partial M\times [0,1)$. By the definition of $h(x_n)\in C^{\infty}([0,1))$
and $h(x_n)>0$, there exists $\tilde{h}\in C^{\infty}((-\varepsilon,1))$ such that $\tilde{h}|_{[0,1)}=h$ and $\tilde{h}>0$ for some
sufficiently small $\varepsilon>0$. Then there exists a metric $\widehat{g}$ on $\widehat{M}=M\bigcup_{\partial M}\partial M\times
(-\varepsilon,0]$ which has the form on $U\bigcup_{\partial M}\partial M\times (-\varepsilon,0 ]$
$$
\widehat{g}=\frac{1}{\tilde{h}(x_{n})}g^{\partial M}+\texttt{d}x _{n}^{2} ,\eqno(2.2)$$
such that $\widehat{g}|_{M}=g$.
We fix a metric $\widehat{g}$ on the
$\widehat{M}$ such that $\widehat{g}|_M=g$. We can get the spin
structure on $\widehat{M}$ by extending the spin structure on $M.$
Let $D$ be the Dirac operator associated to $\widehat{g}$ on the
spinors bundle $S(T\widehat{M})$. We want to compute
 $$\widetilde{{\rm Wres}}\big[\pi^+D^{-1}\circ \pi^+D^{(-2m+1)}\big]$$ (for the related definitions, see , Section 2, 3 [Wa1]).

To define the lower dimensional volume, some basic facts and formulae about Boutet de Monvel's calculus which can be found in Sec.2 in [Wa1]
are needed.
 Let  $$
  F:L^2({\bf R}_t)\rightarrow L^2({\bf R}_v);~F(u)(v)=\int e^{-ivt}u(t)dt
  $$
   denote the Fourier transformation and
$\Phi(\overline{{\bf R}^+}) =r^+\Phi({\bf R})$ (similarly define $\Phi(\overline{{\bf R}^-}$)), where $\Phi({\bf R})$
denotes the Schwartz space and
 $$
r^{+}:C^\infty ({\bf R})\rightarrow C^\infty (\overline{{\bf R}^+});~ f\rightarrow f|\overline{{\bf R}^+};~
 \overline{{\bf R}^+}=\{x\geq0;x\in {\bf R}\}.
$$
We define $H^+=F(\Phi(\overline{{\bf R}^+}));~ H^-_0=F(\Phi(\overline{{\bf R}^-}))$ which are orthogonal to each other. We have the following
 property: $h\in H^+~(H^-_0)$ iff $h\in C^\infty({\bf R})$ which has an analytic extension to the lower (upper) complex
half-plane $\{{\rm Im}\xi<0\}~(\{{\rm Im}\xi>0\})$ such that for all nonnegative integer $l$,
 $$
\frac{d^{l}h}{d\xi^l}(\xi)\sim\sum^{\infty}_{k=1}\frac{d^l}{d\xi^l}(\frac{c_k}{\xi^k})
$$
as $|\xi|\rightarrow +\infty,{\rm Im}\xi\leq0~({\rm Im}\xi\geq0)$.

 Let $H'$ be the space of all polynomials and $H^-=H^-_0\bigoplus H';~H=H^+\bigoplus H^-.$ Denote by $\pi^+~(\pi^-)$ respectively the
 projection on $H^+~(H^-)$. For calculations, we take $H=\widetilde H=\{$rational functions having no poles on the real axis$\}$ ($\tilde{H}$
 is a dense set in the topology of $H$). Then on $\tilde{H}$,
$$
\pi^+h(\xi_0)=\frac{1}{2\pi i}\lim_{u\rightarrow 0^{-}}\int_{\Gamma^+}\frac{h(\xi)}{\xi_0+iu-\xi}d\xi,\eqno(2.3)
$$
where $\Gamma^+$ is a Jordan close curve
included ${\rm Im}(\xi)>0$ surrounding all the singularities of $h$ in the upper half-plane and
$\xi_0\in {\bf R}$. Similarly, define $\pi^{'}$ on $\tilde{H}$,
$$
\pi'h=\frac{1}{2\pi}\int_{\Gamma^+}h(\xi)d\xi.\eqno(2.4)
$$
So, $\pi'(H^-)=0$. For $h\in H\bigcap L^1(R)$, $\pi'h=\frac{1}{2\pi}\int_{R}h(v)dv$ and for $h\in H^+\bigcap L^1(R)$, $\pi'h=0$.

Let $M$ be an $n$-dimensional compact oriented manifold with boundary $\partial M$.
Denote by $\mathcal{B}$ Boutet de Monvel's algebra, we recall the main theorem in [FGLS].\\

\noindent{\bf Theorem 2.1}{\bf(Fedosov-Golse-Leichtnam-Schrohe)}
 {\it Let $X$ and $\partial X$ be connected, ${\rm dim}X=n\geq3$,
 $A=\left(\begin{array}{lcr}\pi^+P+G &   K \\
T &  S    \end{array}\right)$ $\in \mathcal{B}$ , and denote by $p$, $b$ and $s$ the local symbols of $P,G$ and $S$ respectively.
 Define:}
\begin{eqnarray*}
{\rm{\widetilde{Wres}}}(A)&=&\int_X\int_{\bf S}{\rm{tr}}_E\left[p_{-n}(x,\xi)\right]\sigma(\xi)dx \nonumber\\
&&+2\pi\int_ {\partial X}\int_{\bf S'}\left\{{\rm tr}_E\left[({\rm{tr}}b_{-n})(x',\xi')\right]+{\rm{tr}}
_F\left[s_{1-n}(x',\xi')\right]\right\}\sigma(\xi')dx',~~~  (2.5)
\end{eqnarray*}
{\it Then~~ a) ${\rm \widetilde{Wres}}([A,B])=0 $, for any
$A,B\in\mathcal{B}$;~~ b) It is a unique continuous trace on
$\mathcal{B}/\mathcal{B}^{-\infty}$.}\\

\indent Let $p_{1},p_{2}$ be nonnegative integers and $p_{1}+p_{2}\leq n$. Then by Sec 2.1 of [Wa2],  we have\\

\noindent {\bf Definition 2.2.} Lower-dimensional volumes of spin manifolds with boundary  are defined by
  $$\label{}
  {\rm Vol}^{(p_1,p_2)}_nM:=\widetilde{{\rm Wres}}[\pi^+D^{-p_1}\circ\pi^+D^{-p_2}].\eqno(2.6)
$$
Denote by $\sigma_{l}(A)$ the $l$-order symbol of an operator A. Similar to (2.1.7) in [Wa2], we have that
$$
\widetilde{{\rm Wres}}[\pi^+D^{-p_1}\circ\pi^+D^{-p_2}]=\int_M\int_{|\xi|=1}{\rm
trace}_{S(TM)}[\sigma_{-n}(D^{-p_1-p_2})]\sigma(\xi)\texttt{d}x+\int_{\partial
M}\Phi,\eqno(2.7)
$$
where
 \begin{eqnarray}
\Phi&=&\int_{|\xi'|=1}\int^{+\infty}_{-\infty}\sum^{\infty}_{j, k=0}
\sum\frac{(-i)^{|\alpha|+j+k+1}}{\alpha!(j+k+1)!}
 {\rm trace}_{S(TM)}
\Big[\partial^j_{x_n}\partial^\alpha_{\xi'}\partial^k_{\xi_n}
\sigma^+_{r}(D^{-p_1})(x',0,\xi',\xi_n)\nonumber\\
&&\times\partial^\alpha_{x'}\partial^{j+1}_{\xi_n}\partial^k_{x_n}\sigma_{l}
(D^{-p_2})(x',0,\xi',\xi_n)\Big]d\xi_n\sigma(\xi')\texttt{d}x',~~~~~~~~~~~~~~~~~~~~~~~~~~~~~~~~~(2.8)\nonumber
\end{eqnarray}
and the sum is taken over $r-k+|\alpha|+\ell-j-1=-n,r\leq-p_{1},\ell\leq-p_{2}$.
Since $[\sigma_{-n}(D^{-n+2})]|_M$
has the same expression as $\sigma_{-n}(D^{-n+2})$ in the case of
manifolds without boundary, by (2.5) in [Ac], we have
 $$
{\rm Wres}(D^{-n+2})
= \frac{(2-n)(2\pi)^{{n}/2}}{12\Gamma({n}/2)}\int_Ms{\rm dvol}_M.\eqno(2.9)$$
So we only need to compute $\int_{\partial M}\Phi$.

 \section{ Kastler-Kalau-Walze type theorems for any dimensional manifolds with  boundary }
\quad In this section, we compute the lower dimensional volume for any dimensional compact manifolds with boundary and get a
Kastler-Kalau-Walze type formula in this case. From now on we always assume that $M$ carries a spin structure so that the spinor bundle and the
Dirac operator are defined on $M$.

 Set  $\widetilde{E}_{n}=\frac{\partial}{\partial x_{n}}$,
 $\widetilde{E}_{j}=\sqrt{h(x_{n})}E_{j}~~(1\leq j \leq n-1)$, where  $\{E_{1},\cdots,E_{n-1}\}$ are orthonormal basis of $T\partial_{M}$.
  Let $\nabla^L$ denote the Levi-civita connection
about $g^M$.
 In the local coordinates $\{x_i; 1\leq i\leq n\}$ and the fixed orthonormal frame $\{\widetilde{E}_{1},\cdots,\widetilde{E}_{n}\}$,
 the connection matrix
 $(\omega_{s,t})$
is defined by
 $$
 \nabla^L(\widetilde{E}_{1},\cdots,\widetilde{E}_{n})^{t}= (\omega_{s,t})(\widetilde{E}_{1},\cdots,\widetilde{E}_{n})^{t}.$$
 The Dirac operator is defined by
$$
D=\sum^n_{j=1}c(\widetilde{E_{j}})\Big[\widetilde{E_j}+\frac{1}{4}\sum_{s,t}\omega_{s,t}(\widetilde{E_j})c(\widetilde{E_s})c(\widetilde{E_t})\Big].\eqno(3.1)
$$

Let $g^{ij}=g(dx_i,dx_j)$ and
$$
\nabla^L_{\partial_i}\partial_j=\sum_k\Gamma_{ij}^k\partial_k;
~\Gamma^k=g^{ij}\Gamma_{ij}^k, 
 $$
Let the cotangent vector $\xi=\sum \xi_jdx_j$ and
$\xi^j=g^{ij}\xi_i$.
Let the symbols of $D^2$ and $D$ be
 $$
\sigma(D^{2})=\sigma_{2}(D^{2})+\sigma_{1}(D^{2})+\sigma_{0}(D^{2});~~ \sigma(D^{1})=\sigma_{1}(D^{1})+\sigma_{0}(D^{1}).\eqno(3.2)$$
By the composition formula of psudodifferential operators, then we have
$$
\sigma_{-2}(D^{-2})=|\xi|^{-2},~~\sigma_{-(2m-2)}(D^{-(2m-2)})=(|\xi|^2)^{1-m},~~\sigma_{-2m+1}(D^{-2m+1})=\frac{\sqrt{-1}c(\xi)}{|\xi|^{2m}}.\eqno(3.3)$$
By (3.8) in [WW2], we have:

$$\sigma_{-1-2m}(D^{-2m})=m\sigma_2(D^2)^{(-m+1)}\sigma_{-3}(D^{-2})$$
$$-\sqrt{-1} \sum_{k=0}^{m-2}\sum_{\mu=1}^{2m+2}
\partial_{\xi_{\mu}}\sigma_{2}^{-m+k+1}(D^2)
\partial_{x_{\mu}}\sigma_{2}^{-1}(D^2)(\sigma_2(D^2))^{-k}.\eqno(3.4)$$
By (3.4), we know that
$$\sigma_{-2m}(D^{1-2m})=\sigma_{-2m}(D^{-2m}\cdot D)=\left\{\sum_{|\alpha|=0}^{+\infty}(-\sqrt{-1})^{|\alpha|}\frac{1}{\alpha!}\partial^\alpha_\xi[\sigma(D^{-2m})]\partial^\alpha_x[\sigma(D)]\right\}_{-2m}$$
$$=\sigma_{-2m}(D^{-2m})\sigma_0(D)+\sigma_{-2m-1}(D^{-2m})\sigma_1(D)+\sum_{|\alpha|=1}(-\sqrt{-1})
\partial^\alpha_\xi[\sigma_{-2m}(D^{-2m})]\partial^\alpha_x[\sigma_1(D)]$$
$$=|\xi|^{-2m}\sigma_0(D)+\sum_{j=1}^{2m+2}\partial_{\xi_j}(|\xi|^{-2m})\partial_{x_j}c(\xi)+
\left[m\sigma_2(D^2)^{(-m+1)}\sigma_{-3}(D^{-2})\right.$$
$$\left.-\sqrt{-1} \sum_{k=0}^{m-2}\sum_{\mu=1}^{2m+2}
\partial_{\xi_{\mu}}\sigma_{2}^{-m+k+1}(D^2)
\partial_{x_{\mu}}\sigma_{2}^{-1}(D^2)(\sigma_2(D^2))^{-k}\right]\sqrt{-1}c(\xi).
\eqno(3.5)$$
Since $\Phi$ is a global form on $\partial M$, so for any fixed point $x_{0}\in\partial M$, we can choose the normal coordinates
$U$ of $x_{0}$ in $\partial M$(not in $M$) and compute $\Phi(x_{0})$ in the coordinates $\widetilde{U}=U\times [0,1)$ and the metric
$\frac{1}{h(x_{n})}g^{\partial M}+\texttt{d}x _{n}^{2}$. The dual metric of $g^{M}$ on $\widetilde{U}$ is
$h(x_{n})g^{\partial M}+\texttt{d}x _{n}^{2}.$  Now we can compute $\Phi$ (see formula (2.8) for the definition of $\Phi$).
Since the sum is taken over $r+\ell-k-j-|\alpha|-1=-n, \ r\leq-1, \ell\leq 1-2m$, then we have the $\int_{\partial_{M}}\Phi$
is the sum of the following five cases:\\

\noindent  {\bf case a)~I)}~$r=-1,~l=1-2m,~k=j=0,~|\alpha|=1$\\
$${\rm case~a)~I)}=-\int_{|\xi'|=1}\int^{+\infty}_{-\infty}\sum_{|\alpha|=1}{\rm trace}
\Big[\partial^\alpha_{\xi'}\pi^+_{\xi_n}\sigma_{-1}(D^{-1})$$
$$\times\partial^\alpha_{x'}
\partial_{\xi_n}\sigma_{1-2m}(D^{1-2m})\Big](x_0)d\xi_n\sigma(\xi')dx'.\eqno(3.6)$$
By Lemma 2.2 in [Wa2] and (3.12) in [WW2], we have for $j<n$
$$\partial_{x_j}\sigma_{1-2m}(D^{1-2m})(x_0)=\partial_{x_j}[\sqrt{-1}c(\xi)|\xi|^{-2m}]$$
$$=\sqrt{-1}[\partial_{x_j}c(\xi)](x_0)|\xi|^{-2m}+\sqrt{-1}c(\xi)\partial_{x_j}(|\xi|^{-2m})(x_0)=0.\eqno(3.7)$$
 so {\rm {\bf case~a)~I)}} vanishes.\\

\noindent  {\bf case a)~II)}~$r=-1,~l=1-2m,~k=|\alpha|=0,~j=1$\\

\indent By (2.8), we get
 $$
{\rm case~a)~II)}=-\frac{1}{2}\int_{|\xi'|=1}\int^{+\infty}_{-\infty} {\rm
trace} \Big[\partial_{x_n}\pi^+_{\xi_n}\sigma_{-1}(D^{-1})$$
$$\times \partial_{\xi_n}^2
\sigma_{1-2m}(D^{1-2m})\Big](x_0)d\xi_n\sigma(\xi')dx'$$
$$=-\frac{1}{2}\int_{|\xi'|=1}\int^{+\infty}_{-\infty} {\rm
trace} \Big[ \partial_{\xi_n}^2\partial_{x_n}\pi^+_{\xi_n}\sigma_{-1}(D^{-1})$$
$$\times
\sigma_{1-2m}(D^{1-2m})\Big](x_0)d\xi_n\sigma(\xi')dx'
.\eqno(3.8)
$$
By (2.2.23) in [Wa2], we have
$$\pi^+_{\xi_n}\partial_{x_n}\sigma_{-1}(D^{-1})(x_0)|_{|\xi'|=1}=\frac{\partial_{x_n}[c(\xi')](x_0)}{2(\xi_n-i)}+\sqrt{-1}h'(0)
\left[\frac{ic(\xi')}{4(\xi_n-i)}+\frac{c(\xi')+ic(dx_n)}{4(\xi_n-i)^2}\right].\eqno(3.9)$$
So
$$\partial_{\xi_n}^2\pi^+_{\xi_n}\partial_{x_n}\sigma_{-1}(D^{-1})(x_0)|_{|\xi'|=1}=\frac{\partial_{x_n}[c(\xi')](x_0)}{4(\xi_n-i)^3}+\sqrt{-1}h'(0)
\left[\frac{ic(\xi')}{8(\xi_n-i)^3}+\frac{c(\xi')+ic(dx_n)}{24(\xi_n-i)^4}\right].\eqno(3.10)$$
\noindent We know that
$$\sigma_{1-2m}(D^{1-2m})=\frac{\sqrt{-1}[c(\xi')+\xi_nc(dx_n)]}{(1+\xi_n^2)^m},\eqno(3.11)$$
\noindent By the relation of the Clifford action and ${\rm tr}{AB}={\rm tr }{BA}$, then we have the equalities:\\
$${\rm tr}[c(\xi')c(dx_n)]=0;~~{\rm tr}[c(dx_n)^2]=-2^{m+1};~~{\rm tr}[c(\xi')^2](x_0)|_{|\xi'|=1}=-2^{m+1};~~$$
$${\rm tr}[\partial_{x_n}c(\xi')c(dx_n)]=0;~~{\rm tr}[\partial_{x_n}c(\xi')c(\xi')](x_0)|_{|\xi'|=1}=-2^mh'(0).\eqno(3.12)$$
By (3.10),(3.11) and (3.12), we have
\begin{eqnarray*}
&& {\rm
trace} \Big[ \partial_{\xi_n}^2\partial_{x_n}\pi^+_{\xi_n}\sigma_{-1}(D^{-1})
\times
\sigma_{1-2m}(D^{1-2m})\Big](x_0)|_{|\xi'|=1}\\
&=& {\rm
trace}\left\{\left\{\frac{\partial_{x_n}[c(\xi')](x_0)}{4(\xi_n-i)^3}+\sqrt{-1}h'(0)
\left[\frac{ic(\xi')}{8(\xi_n-i)^3}+\frac{c(\xi')+ic(dx_n)}{24(\xi_n-i)^4}\right]\right\}\right.\\
&&\left.\times
\frac{\sqrt{-1}[c(\xi')+\xi_nc(dx_n)]}{(1+\xi_n^2)^m}\right\}\\
&=& {\rm
trace}\left\{\left[\frac{\partial_{x_n}[c(\xi')](x_0)}{4(\xi_n-i)^3}+
\frac{(4i-3\xi_n)h'(0)}{24(\xi_n-i)^4}c(\xi')-\frac{h'(0)c(dx_n)}{24(\xi_n-i)^4}\right]\right.\\
&&\left.\times
\frac{\sqrt{-1}[c(\xi')+\xi_nc(dx_n)]}{(1+\xi_n^2)^m}\right\}\\
&=&\frac{2^mh'(0)i}{12(\xi_n+i)^m(\xi_n-i)^{m+3}}.~~~~~~~~~~~~~~~~~~~~~~~~~~~~~~~~~~~~~~~~~~~~~~~~~~~~~~~~~~~~~~~~ (3.13)
\end{eqnarray*}
By (3.13) and the Cauchy integral formula, we have
$$
{\rm {\bf case~a)~II)}}=Vol(S_{n-2}) \frac{2^mh'(0)\pi}{12(m+2)!}[(\xi_n+i)^{-m}]^{(m+2)}|_{\xi_n=i}dx',\eqno(3.14)
$$
 where $Vol(S_{n-2})$ is the canonical volume of $S_{n-2}$ and denote the $p$-th derivative of $f(\xi_n)$ by $[f(\xi_n)]^{(p)}$.\\

\noindent  {\bf case a)~III)}~$r=-1,~l=1-2m,~j=|\alpha|=0,~k=1$\\

 \indent By (2.8) and an integration by parts, we get
 \begin{eqnarray}
{\rm case~ a)~III)}&=&-\frac{1}{2}\int_{|\xi'|=1}\int^{+\infty}_{-\infty}{\rm trace} \Big[\partial_{\xi_n}\pi^+_{\xi_n}\sigma_{-1}(D^{-1})\nonumber\\
&&\times \partial_{\xi_n}\partial_{x_n}\sigma_{1-2m}(D^{1-2m})\Big](x_0)d\xi_n\sigma(\xi')dx'\nonumber\\
                   &=& \frac{1}{2}\int_{|\xi'|=1}\int^{+\infty}_{-\infty} {\rm trace}\Big[\partial_{\xi_n}^2\pi^+_{\xi_n}\sigma_{-1}(D^{-1})\nonumber\\
                   &&\times\partial_{x_n}\sigma_{1-2m}(D^{1-2m})\Big](x_0)d\xi_n\sigma(\xi')dx'.~~~~~~~~~~~~~~~~~~~~~~~~~~~~~~~~~ (3.15)\nonumber
\end{eqnarray}
By (2.2.29) in [Wa2], we have
$$\partial_{\xi_n}^2\pi^+_{\xi_n}\sigma_{-1}(D^{-1})(x_0)|_{|\xi'|=1}=\frac{c(\xi')+ic(dx_n)}{(\xi_n-i)^3}.\eqno(3.16)$$
By (3.3), direct computations show that
$$\partial_{x_n}\sigma_{1-2m}(D^{1-2m})(x_0)|_{|\xi'|=1}=\frac{\sqrt{-1}\partial_{x_n}[c(\xi')](x_0)}
{(1+\xi_n^2)^m}-\frac{\sqrt{-1}mh'(0)c(\xi)}{(1+\xi^2_n)^{m+1}}.\eqno(3.17)$$
So by (3.16),(3.17) and (3.12) and the Cauchy integral formula, we have
$${\rm case~ a)~III)}=\frac{1}{2}\int_{|\xi'|=1}\int^{+\infty}_{-\infty} {\rm trace}\left\{\frac{c(\xi')+ic(dx_n)}{(\xi_n-i)^3}\right.$$
$$\left.\times
\left[\frac{\sqrt{-1}\partial_{x_n}[c(\xi')](x_0)}
{(1+\xi_n^2)^m}-\frac{\sqrt{-1}mh'(0)c(\xi)}{(1+\xi^2_n)^{m+1}}\right]\right\}(x_0)d\xi_n\sigma(\xi')dx'$$
$$=\frac{1}{2}\int_{|\xi'|=1}\int^{+\infty}_{-\infty}2^mh'(0)\times\frac{-i\xi^2_n-2m\xi_n+2mi-i}{(\xi_n+i)^{m+1}(\xi_n-i)^{m+4}}d\xi_n\sigma(\xi')dx'$$
$$=\frac{\pi ih'(0)2^m{\rm Vol}(S^{n-2})dx'}{(m+3)!}\left[\frac{-i\xi^2_n-2m\xi_n+2mi-i}{(\xi_n+i)^{m+1}}\right]^{(m+3)}|_{\xi_n=i}.\eqno(3.18)$$\\

\noindent  {\bf case b)}~$r=-2,~l=1-2m,~k=j=|\alpha|=0$\\

\indent By (2.8), we get
$$
{\rm case~ b)}=-i\int_{|\xi'|=1}\int^{+\infty}_{-\infty}{\rm trace} [\pi^+_{\xi_n}\sigma_{-2}(D^{-1})\times
     \partial_{\xi_n}\sigma_{1-2m}(D^{1-2m})](x_0)d\xi_n\sigma(\xi')dx'.\eqno(3.19)$$
\indent By (2.2.34)-(2.2.37) in [Wa2], we have
$$\pi^+_{\xi_n}\sigma_{-2}(D^{-1})(x_0)|_{|\xi'|=1}
=B_1-B_2,\eqno (3.20)$$
\noindent where\\
$$B_1=-\frac{A_1}{4(\xi_n-i)}-\frac{A_2}{4(\xi_n-i)^2},\eqno(3.21)$$
\noindent and
$$A_1=ic(\xi')\sigma_0(D)c(\xi')+ic(dx_n)\sigma_0(D)c(dx_n)+ic(\xi')c(dx_n)\partial_{x_n}[c(\xi')];$$
$$A_2=[c(\xi')+ic(dx_n)]\sigma_0(D)[c(\xi')+ic(dx_n)]+c(\xi')c(dx_n)\partial_{x_n}c(\xi')-i\partial_{x_n}[c(\xi')].\eqno(3.22)$$
$$
B_2
=\frac{h'(0)}{2}\left[\frac{c(dx_n)}{4i(\xi_n-i)}+\frac{c(dx_n)-ic(\xi')}{8(\xi_n-i)^2}
+\frac{3\xi_n-7i}{8(\xi_n-i)^3}[ic(\xi')-c(dx_n)]\right].\eqno(3.23)$$
Similar to (2.2.38) in [Wa2], we have
$$\partial_{\xi_n}\sigma_{1-2m}(D^{1-2m})(x_0)|_{|\xi'|=1}=\sqrt{-1}\left[\frac{c(dx_n)}{(1+\xi_n^2)^m}-m\times\frac{2\xi_nc(\xi')+2\xi_n^2c(dx_n)}
{(1+\xi_n^2)^{m+1}}\right].\eqno(3.24)$$
\noindent By (3.23), (3.24) and (3.12), we have\\
$${\rm tr }[B_2\times\partial_{\xi_n}\sigma_{1-2m}(D^{1-2m})(x_0)]|_{|\xi'|=1}
=\frac{\sqrt{-1}}{2}h'(0){\rm trace}$$
$$\left\{\{\left[\frac{1}{4i(\xi_n-i)}+\frac{1}{8(\xi_n-i)^2}-\frac{3\xi_n-7i}{8(\xi_n-i)^3}\right]c(dx_n)
+\left[\frac{-1}{8(\xi_n-i)^2}+\frac{3\xi_n-7i}{8(\xi_n-i)^3}\right]ic(\xi')\}\right.$$
$$\times\left.\{\left[\frac{1}{(1+\xi_n^2)^{m}}-\frac{2m\xi_n^2}{(1+\xi_n^2)^{m+1}}\right]c(dx_n)-\frac{2m\xi_n}
{(1+\xi_n^2)^{m+1}}c(\xi')\}\right\}$$
$$=h'(0)2^{m-2}\times \frac{(2m-1)\xi^3_n-2i(2m-1)\xi^2_n-(6m-1)\xi_n+4i}{(\xi_n-i)^2(1+\xi^2_n)^{m+1}}.\eqno(3.25)$$

\noindent By (2.2.40) in [Wa2], we have
$$B_1=\frac{-1}{4(\xi_n-i)^2}[(2+i\xi_n)c(\xi')\sigma_0(D)c(\xi')+i\xi_nc(dx_n)\sigma_0(D)c(dx_n)$$
$$+
(2+i\xi_n)c(\xi')c(dx_n)\partial_{x_n}c(\xi')+ic(dx_n)\sigma_0(D)c(\xi')
+ic(\xi')\sigma_0(D)c(dx_n)-i\partial_{x_n}c(\xi')].\eqno(3.26)$$
By (3.24), we have
$$\partial_{\xi_n}\sigma_{1-2m}(D^{1-2m})(x_0)|_{|\xi'|=1}=\sqrt{-1}\left[\frac{1+(1-2m)\xi^2_n}{(1+\xi_n^2)^{m+1}}c(dx_n)
-\frac{2m\xi_n}
{(1+\xi_n^2)^{m+1}}c(\xi')\right].\eqno(3.27)$$
Similar to Lemma 2.4 in [Wa2], we have
$$\sigma_0(D)(x_0)=c_0c(dx_n),~~{\rm where}~ c_0=\frac{1-n}{4}h'(0).\eqno(3.28)$$
\noindent By the relation of the Clifford action and ${\rm tr}{AB}={\rm tr }{BA}$, then we have the equalities:\\
$${\rm tr}[c(\xi')\sigma_0(D)c(\xi')c(dx_n)]=-c_02^{m+1};~~{\rm tr}[c(dx_n)\sigma_0(D)c(dx_n)^2]=c_02^{m+1};$$
$${\rm tr}[c(\xi')c(dx_n)\partial_{x_n}c(\xi')c(dx_n)](x_0)|_{|\xi'|=1}=-2^{m}h'(0);~
{\rm tr}[c(dx_n)\sigma_0(D)c(\xi')^2]=c_02^{m+1}.\eqno(3.29)$$

\noindent By (3.26)-(3.29),
 considering for $i<n$, $\int_{|\xi'|=1}\{{\rm odd ~number~ product~ of~}\xi_i\}\sigma(\xi')=0$, then\\
$${\rm tr }[B_1\times\partial_{\xi_n}\sigma_{1-2m}(D^{1-2m})(x_0)]|_{|\xi'|=1}=\frac{2^mih'(0)}{4(\xi_n-i)^2(1+\xi_n^2)^{m+1}}$$
$$\cdot\left\{
(n-1)[(2m-1)\xi_n^2-2im\xi_n-1]\right.$$
$$\left.+[(1-2m)i\xi_n^3+2(1-2m)\xi_n^2+(2m+1)i\xi_n+2]\right\}.\eqno(3.30)$$
By (3.25) and (3.30) and the Cauchy integral formula, we have
$$
{\rm case~ b)}=-i\int_{|\xi'|=1}\int^{+\infty}_{-\infty}{\rm trace} [(B_1-B_2)\times
     \partial_{\xi_n}\sigma_{1-2m}(D^{1-2m})](x_0)d\xi_n\sigma(\xi')dx'$$
 $$=2^mh'(0)\int_{|\xi'|=1}\int^{+\infty}_{-\infty}\frac{(4m^2-1)\xi_n^2+(-4m^2-6m+2)i\xi_n-(n+1)}{4(\xi_n-i)^{m+3}(\xi_n+i)^{m+1}}
     d\xi_n\sigma(\xi')dx'$$
$$=\frac{2^{m+1}h'(0){\rm Vol}(S^{n-2})\pi i dx'}{(m+2)!}\left[
    \frac{(4m^2-1)\xi_n^2+(-4m^2-6m+2)i\xi_n-(n+1)}{4(\xi_n+i)^{m+1}}\right]^{(m+2)}|_{\xi_n=i}.\eqno(3.31)$$ \\

  \noindent  {\bf case c)}~$r=-1,~l=-2m,~k=j=|\alpha|=0$\\

\indent By (2.8) and an integration by parts, we get
$$
{\rm case~ c)}=-i\int_{|\xi'|=1}\int^{+\infty}_{-\infty}{\rm trace} [\pi^+_{\xi_n}\sigma_{-1}(D^{-1})\times
     \partial_{\xi_n}\sigma_{-2m}(D^{-2m+1})](x_0)d\xi_n\sigma(\xi')dx'$$
    $$ =i\int_{|\xi'|=1}\int^{+\infty}_{-\infty}{\rm trace} [\partial_{\xi_n}\pi^+_{\xi_n}\sigma_{-1}(D^{-1})\times
        \sigma_{-2m}(D^{-2m+1})](x_0)d\xi_n\sigma(\xi')dx'.\eqno(3.32)
$$
By (2.2.44) in [Wa2] and similar to (3.29) in [WW2], we have
$$ \partial_{\xi_n}\pi^+_{\xi_n}\sigma_{-1}(D^{-1})(x_0)|_{|\xi'|=1}=-\frac{c(\xi')+ic(dx_n)}{2(\xi_n-i)^2};\eqno(3.33)$$
By Lemma 2.2 in [Wa2], we have
$$\sum_{j=1}^{2m+2}\partial_{\xi_j}(|\xi|^{-2m})\partial_{x_j}(c(\xi))(x_0)|_{|\xi'|=1}=-2m\xi_n(1+\xi^2_n)^{-m-1}\partial_{x_n}[c(\xi')](x_0).\eqno(3.34)$$
$$\left[-\sqrt{-1} \sum_{k=0}^{m-2}\sum_{\mu=1}^{2m+2}
\partial_{\xi_{\mu}}\sigma_{2}^{-m+k+1}(D^2)
\partial_{x_{\mu}}\sigma_{2}^{-1}(D^2)(\sigma_2(D^2))^{-k}\right]\sqrt{-1}c(\xi)(x_0)|_{|\xi'|=1}$$
$$
=-2c(\xi)h'(0)\xi_n(-\frac{m^2}{2}+\frac{m}{2})(1+\xi_n^2)^{-m-2}.\eqno(3.35)$$
By (3.26) in [WW2], we have
$$
\sigma_{-3}(D^{-2})(x_0)|_{|\xi'|=1}
                                   =\frac{-i}{(1+\xi_n^2)^2}\Big(-\frac{1}{2}h'(0)\sum_{k<n}\xi_k c(\widetilde{e_k})c(\widetilde{e_n})
                                       +\frac{n-1}{2}h'(0)\xi_n\Big)-\frac{2ih'(0)\xi_n}{(1+\xi_n^2)^3}$$
$$=\frac{-i}{(1+\xi_n^2)^2}\Big(-\frac{1}{2}h'(0)c(\xi')c(dx_n)
                                       +\frac{n-1}{2}h'(0)\xi_n\Big)-\frac{2ih'(0)\xi_n}{(1+\xi_n^2)^3}.\eqno(3.36)$$
So by (3.5), (3.34)-(3.36), we have
$$\sigma_{-2m}(D^{1-2m})(x_0)|_{|\xi'|=1}=\frac{(1-n)h'(0)c(dx_n)}{4(1+\xi_n^2)^m}-2m\xi_n(1+\xi^2_n)^{-m-1}\partial_{x_n}[c(\xi')](x_0)$$
$$
+mi(1+\xi_n^2)^{-m+1}[c(\xi')+\xi_nc(dx_n)]\times\left[
\frac{-ih'(0)c(\xi')c(dx_n)}{2(1+\xi_n^2)^2}
                                       +\frac{(-n+1)h'(0)i\xi_n}{2(1+\xi_n^2)^2}-\frac{2ih'(0)\xi_n}{(1+\xi_n^2)^3}\right]$$
$$-[c(\xi')+\xi_nc(dx_n)]h'(0)\xi_n(-{m^2}+{m})(1+\xi_n^2)^{-m-2}.\eqno(3.37)$$
By (3.33) and (3.37), we have
$${\rm trace} [\partial_{\xi_n}\pi^+_{\xi_n}\sigma_{-1}(D^{-1})\times
        \sigma_{-2m}(D^{-2m+1})](x_0)|_{|\xi'|=1}=\frac{2^mh'(0)}
        {4(\xi_n-i)^{m+3}(\xi_n+i)^{m+2}}$$
        $$\times
        \left[(2nm-2m-n+1)i\xi_n^3+(-2m+1+2\pi)\xi_n^2\right.$$
        $$\left.+(2nm-2m+4m^2-n+1)i\xi_n+(n-1+2\pi)\right]
        .\eqno(3.38)$$
Then by the Cauchy integral formula, we get
$$
{\rm case~ c)}=\frac{-\pi 2^{m-1}h'(0){\rm Vol }(S^{n-2})dx'}{(m+2)!}\left\{\frac{1}{(\xi_n+i)^{m+2}}
\left[(2nm-2m-n+1)i\xi_n^3\right.\right.$$
$$\left.\left.+(-2m+1+2\pi)\xi_n^2
        +(2nm-2m+4m^2-n+1)i\xi_n+(n-1+2\pi)\right]\right\}|_{|\xi_n=i}^{(m+2)}.\eqno(3.39)$$
Similar to (3.41) in [WW2], we have for the extrinsic curvature
$$K(x_0)=\frac{1-n}{2}h'(0).\eqno(3.40)$$
By (3.14),(3.18),(3.31), (3.39) and (3.40), we have
$$\Phi=\frac{2}{1-n}{\rm Vol}(S^{n-2})L_0\int_{\partial M}K{\rm dvol}_{\partial M},\eqno(3.41)$$
where $$L_0=
 \frac{2^m\pi}{12(m+2)!}[(\xi_n+i)^{-m}]^{(m+2)}|_{\xi_n=i}
+\frac{\pi i2^m}{(m+3)!}\left[\frac{-i\xi^2_n-2m\xi_n+2mi-i}{(\xi_n+i)^{m+1}}\right]^{(m+3)}|_{\xi_n=i}$$
$$
+\frac{2^{m+1}\pi i }{(m+2)!}\left[
    \frac{(4m^2-1)\xi_n^2+(-4m^2-6m+2)i\xi_n-(n+1)}{4(\xi_n+i)^{m+1}}\right]^{(m+2)}|_{\xi_n=i}$$
    $$
+\frac{-\pi 2^{m-1}}{(m+2)!}\left\{\frac{1}{(\xi_n+i)^{m+2}}
\left[(2nm-2m-n+1)i\xi_n^3\right.\right.$$
$$\left.\left.+(-2m+1+2\pi)\xi_n^2
        +(2nm-2m+4m^2-n+1)i\xi_n+(n-1+2\pi)\right]\right\}|_{|\xi_n=i}^{(m+2)}.\eqno(3.42)$$
By (2.9), (3.41) and (3.42), we prove Theorem 1.1.\\

\indent Nextly, for $n=2m+1$-dimensional spin manifolds with boundary, we compute
$ \widetilde{{\rm Wres}}\big[\pi^+D^{-1}\circ \pi^+D^{(1-2m)}\big]$.
By Proposition 3.5 in [WW2], we have
 $$
 \widetilde{{\rm Wres}}\big[\pi^+D^{-1}\circ \pi^+D^{(1-2m)}\big]=\int_{\partial M}{\Phi}.\eqno(3.43)
$$
From the formula (2.8) for the definition of $\Phi$, and the sum is taken over
 $r-k+|\alpha|+\ell-j-1=-(2m+1),~~r\leq-1,~~\ell\leq 1-2m$, then
  $r=-1,~~\ell=1-2m,~~k=|\alpha|=j=0$,
$$
{\Phi}=-i\int_{|\xi'|=1}\int^{+\infty}_{-\infty} {\rm trace}_{S(TM)}[\pi^+_{\xi_n} \sigma_{-1}(D^{-1})\times
\partial_{\xi_n}\sigma_{1-2m}(D^{1-2m})]d\xi_n\sigma(\xi')dx'$$
$$=i\int_{|\xi'|=1}\int^{+\infty}_{-\infty} {\rm trace}_{S(TM)}[\partial_{\xi_n}\pi^+_{\xi_n} \sigma_{-1}(D^{-1})\times
\sigma_{1-2m}(D^{1-2m})]d\xi_n\sigma(\xi')dx'
\eqno(3.44)
$$
By (5.3) in [Wa2], we have
$$\partial_{\xi_n}\pi^+_{\xi_n} \sigma_{-1}(D^{-1})(x_0)|_{|\xi'|=1}=-\frac{c(\xi')+ic(dx_n)}{2(\xi_n-i)^2}.\eqno(3.45)$$
So
$$ {\rm trace}_{S(TM)}[\partial_{\xi_n}\pi^+_{\xi_n} \sigma_{-1}(D^{-1})\times
\sigma_{1-2m}(D^{1-2m})]
(x_0)|_{|\xi'|=1}$$
$$=\frac{-i}{2(\xi_n-i)^2(1+\xi^2_n)^m}{\rm trace}_{S(TM)}\left\{[c(\xi')+ic(dx_n)][c(\xi')+\xi_nc(dx_n)]\right\}$$
$$=\frac{-2^{m-1}}{(\xi_n-i)^{m+1}(\xi_n+i)^{m}}.\eqno(3.46)$$
Then
$$\int_{\partial M}\Phi=\frac{2^{-m}m(m+1)\cdots (2m-1){\rm Vol}_{\partial M}{\rm Vol}(S^{n-2})\pi}{m!}.\eqno(3.47)$$
So we have\\

\noindent{\bf Theorem 3.1}
{\it Let M be a $n=2m+1$-dimensional compact spin manifold  with the boundary $\partial M$, then}
 $$\widetilde{{\rm Wres}}\big[\pi^+D^{-1}\circ \pi^+D^{(1-2m)}\big]=\frac{2^{-m}m(m+1)\cdots (2m-1){\rm Vol}_{\partial M}{\rm Vol}(S^{n-2})\pi}{m!}.\eqno(3.48)$$\\

Nextly
for $n=2m+1$-dimensional spin manifolds with boundary, we compute
$ \widetilde{{\rm Wres}}\big[\pi^+D^{-1}\circ \pi^+D^{(2-2m)}\big]$.
 Now we can compute $\Phi$ (see formula (2.8) for the definition of $\Phi$).
Since the sum is taken over $r+\ell-k-j-|\alpha|-1=-(2m+1), \ r\leq-1, \ell\leq 2-2m$, then we have the $\int_{\partial_{M}}\Phi$
is the sum of the following five cases:\\

\noindent  {\bf case a)~I)}~$r=-1,~l=2-2m,~k=j=0,~|\alpha|=1$\\
$${\rm case~a)~I)}=-\int_{|\xi'|=1}\int^{+\infty}_{-\infty}\sum_{|\alpha|=1}{\rm trace}
\Big[\partial^\alpha_{\xi'}\pi^+_{\xi_n}\sigma_{-1}(D^{-1})$$
$$\times\partial^\alpha_{x'}
\partial_{\xi_n}\sigma_{2-2m}(D^{2-2m})\Big](x_0)d\xi_n\sigma(\xi')dx'.\eqno(3.49)$$
By (3.12) in [WW2], we have for $j<n$,
$\partial_{x_j}\sigma_{2-2m}(D^{2-2m})(x_0)=0$. So ${\rm case~a)~I)}=0$.\\

\noindent  {\bf case a)~II)}~$r=-1,~l=2-2m,~k=|\alpha|=0,~j=1$\\

\indent By (2.8), we get
 $$
{\rm case~a)~II)}=-\frac{1}{2}\int_{|\xi'|=1}\int^{+\infty}_{-\infty} {\rm
trace} \Big[\partial_{x_n}\pi^+_{\xi_n}\sigma_{-1}(D^{-1})$$
$$\times \partial_{\xi_n}^2
\sigma_{1-2m}(D^{1-2m})\Big](x_0)d\xi_n\sigma(\xi')dx'
.\eqno(3.50)$$
By (3.9) and (3.16) in [WW2] and $${\rm tr}[c(\xi')]={\rm tr}[c(dx_n)]={\rm tr}[\partial_{x_n}c(\xi')]=0,\eqno(3.51)$$
we have $
{\rm case~a)~II)}=0.$\\

\noindent  {\bf case a)~III)}~$r=-1,~l=2-2m,~j=|\alpha|=0,~k=1$\\

 \indent By (2.8) and an integration by parts, we get
 \begin{eqnarray}
{\rm case~ a)~III)}&=&-\frac{1}{2}\int_{|\xi'|=1}\int^{+\infty}_{-\infty}{\rm trace} \Big[\partial_{\xi_n}\pi^+_{\xi_n}\sigma_{-1}(D^{-1})\nonumber\\
&&\times \partial_{\xi_n}\partial_{x_n}\sigma_{2-2m}(D^{2-2m})\Big](x_0)d\xi_n\sigma(\xi')dx'\nonumber\\
                   &=& \frac{1}{2}\int_{|\xi'|=1}\int^{+\infty}_{-\infty} {\rm trace}\Big[\partial_{\xi_n}^2\pi^+_{\xi_n}\sigma_{-1}(D^{-1})\nonumber\\
                   &&\times\partial_{x_n}\sigma_{2-2m}(D^{2-2m})\Big](x_0)d\xi_n\sigma(\xi')dx'.~~~~~~~~~~~~~~~~~~~~~~~~~~~~~~~~ (3.52)\nonumber
\end{eqnarray}
By (3.16) and (3.21) in [WW2] and (3.51), we have ${\rm case~ a)~III)}=0.$\\

\noindent  {\bf case b)}~$r=-2,~l=2-2m,~k=j=|\alpha|=0$\\

\indent By (2.8), we get
$$
{\rm case~ b)}=-i\int_{|\xi'|=1}\int^{+\infty}_{-\infty}{\rm trace} [\pi^+_{\xi_n}\sigma_{-2}(D^{-1})\times
     \partial_{\xi_n}\sigma_{2-2m}(D^{2-2m})](x_0)d\xi_n\sigma(\xi')dx'.\eqno(3.53)$$
By (3.33) in [WW2] and (3.20)-(3.23), (3.51), (3.53) and
$${\rm tr}[c(\xi')\sigma_0(D)c(\xi')]={\rm tr}[c(dx_n)\sigma_0(D)c(dx_n)]={\rm tr}[c(\xi')c(dx_n)\partial_{x_n}c(\xi')]=0,\eqno(3.54)$$
we have $
{\rm case~ b)}=0.$\\

\noindent  {\bf case c)}~$r=-1,~l=1-2m,~k=j=|\alpha|=0$\\

\indent By (2.8) and an integration by parts, we get
$$
{\rm case~ c)}=-i\int_{|\xi'|=1}\int^{+\infty}_{-\infty}{\rm trace} [\pi^+_{\xi_n}\sigma_{-1}(D^{-1})\times
     \partial_{\xi_n}\sigma_{1-2m}(D^{-2m+2})](x_0)d\xi_n\sigma(\xi')dx'$$
    $$ =i\int_{|\xi'|=1}\int^{+\infty}_{-\infty}{\rm trace} [\partial_{\xi_n}\pi^+_{\xi_n}\sigma_{-1}(D^{-1})\times
        \sigma_{1-2m}(D^{-2m+2})](x_0)d\xi_n\sigma(\xi')dx'.\eqno(3.55)
$$
By (3.33) and (3.8),(3.26),(3.29) in [WW2] and (3.51),(3.54), (3.55), we have $
{\rm case~ c)}=0$. So we get $\int_{\partial_{M}}\Phi=0$.\\

\noindent{\bf Theorem 3.2}
{\it Let M be a $n=2m+1$-dimensional compact spin manifold  with the boundary $\partial M$, then}
$$
 \widetilde{{\rm Wres}}\big[\pi^+D^{-1}\circ \pi^+D^{(-2m+2)}\big]
 = \frac{(2-n)(2\pi)^{{n}/2}}{12\Gamma({n}/2)}\int_Ms{\rm dvol}_M.\eqno(3.56)$$

 \indent Nextly, for $n=2m$-dimensional spin manifolds with boundary, we compute
$ \widetilde{{\rm Wres}}\big[\pi^+D^{-1}\circ \pi^+D^{(2-2m)}\big]$.
By Proposition 3.5 in [WW2], we have
 $$
 \widetilde{{\rm Wres}}\big[\pi^+D^{-1}\circ \pi^+D^{(2-2m)}\big]=\int_{\partial M}{\Phi}.\eqno(3.57)
$$
From the formula (2.8) for the definition of $\Phi$, and the sum is taken over
 $r-k+|\alpha|+\ell-j-1=-2m,~~r\leq-1,~~\ell\leq 2-2m$, then
  $r=-1,~~\ell=2-2m,~~k=|\alpha|=j=0$,
$$
{\Phi}=-i\int_{|\xi'|=1}\int^{+\infty}_{-\infty} {\rm trace}_{S(TM)}[\pi^+_{\xi_n} \sigma_{-1}(D^{-1})\times
\partial_{\xi_n}\sigma_{2-2m}(D^{2-2m})]d\xi_n\sigma(\xi')dx'.
\eqno(3.58)
$$
By (2.2.44) in [Wa2] and (3.47) in [WW2] and (3.51),(3.58), we get $
{\Phi}=0.$\\

\noindent{\bf Theorem 3.3}
{\it Let M be a $n=2m$-dimensional compact spin manifold  with the boundary $\partial M$, then}
$$
 \widetilde{{\rm Wres}}\big[\pi^+D^{-1}\circ \pi^+D^{(-2m+2)}\big]
 =0.\eqno(3.59)$$\\

Similar to Theorem 3.2, we can get\\

\noindent{\bf Theorem 3.4}
{\it Let M be a $n=2m+3$-dimensional compact spin manifold  with the boundary $\partial M$, then}
$$
 \widetilde{{\rm Wres}}\big[\pi^+D^{-2}\circ \pi^+D^{(-2m+1)}\big]
 = \frac{(2-n)(2\pi)^{{n}/2}}{12\Gamma({n}/2)}\int_Ms{\rm dvol}_M.\eqno(3.60)$$

\noindent{\bf References}\\

\noindent [Ac] T. Ackermann, A note on the Wodzicki residue, J.
Geom. Phys., 20, 404-406, 1996.\\
\noindent [Ad] M. Adler.: On a trace functional for formal pseudo-differential operators and the symplectic
structure of Korteweg-de Vries type equations, Invent. Math. 50, 219-248, 1979.\\
\noindent [Co1] A. Connes, Quantized calculus and applications,
XIth International Congress of Mathematical Physics (paris,1994),
15-36, Internat Press, Cambridge, MA, 1995.\\
\noindent [Co2] A. Connes. The action functinal in
noncommutative
geometry, Comm. Math. Phys., 117:673-683, 1998.\\
\noindent [FGLS] B. V. Fedosov, F. Golse, E. Leichtnam, and E.
Schrohe. The noncommutative residue for manifolds with
boundary, J. Funct.
Anal, 142:1-31,1996.\\
 \noindent [Gu] V.W. Guillemin, A new proof of Weyl's
formula on the asymptotic distribution of eigenvalues, Adv. Math.
55 no.2, 131-160, 1985.\\
\noindent [Ka] D. Kastler, The Dirac operator and gravitiation,
Commun. Math. Phys, 166:633-643, 1995.\\
\noindent [KW] W. Kalau and M.Walze, Gravity, non-commutative
geometry, and the Wodzicki residue, J. Geom. Phys., 16:327-344, 1995.\\
\noindent [Po] R. Ponge, Noncommutative Geometry and lower dimensional volumes in Riemannian geometry, Lett. Math. Phys. 83, 1-19, 2008.\\
\noindent [Sc] E. Schrohe, Noncommutative residue, Dixmier's
trace, and
heat trace expansions on manifolds with boundary, Contemp. Math. 242, 161-186, 1999.\\
\noindent [WW1] J. Wang and Y. Wang, The Kastler-Kalau-Walze type theorem for six-dimensional manifolds with boundary,  J. Math. Phys. 56, no. 5, 052501, 14 pp, 2015\\
\noindent [WW2] J. Wang and Y. Wang, A general Kastler-Kalau-Walze type theorem for manifolds with boundary, Int. J. Geom. Methods Mod. Phys. 13, no. 1, 1650003, 16 pp, 2016.\\
\noindent [Wa1] Y. Wang, Diffential forms and the Wodzicki residue for manifolds with boundary. J. Geom. Phys. 56, 731-753, 2006.\\
\noindent [Wa2] Y. Wang, Gravity and the noncommutative residue for manifolds with boundary. Lett. Math. Phys. 80, 37-56, 2007.\\
\noindent [Wa3] Y. Wang, Lower-dimensional volumes and Kastler-kalau-Walze type theorem for manifolds with boundary.
      Commun. Theor. Phys. Vol 54, 38-42, 2010.\\
\noindent [Wo] M. Wodzicki.: Local invariants of spectral asymmetry. Invent. Math. 75(1), 143-178, 1995.\\

\end{document}